\documentclass{Tran-l}
\usepackage{amssymb}
\usepackage{amsmath}
\usepackage{amsfonts}
\usepackage{eurosym}
\usepackage{lscape}
\usepackage[ansinew]{inputenc}
 \newtheorem{thm}{Theorem}[section]
 
 \newtheorem{lem}[thm]{Lemma}
 
 \theoremstyle{definition}
 
 \theoremstyle{remark}
 
 \numberwithin{equation}{section}

\newcommand{\C}{\mathcal{C}}
\newcommand{\D}{\mathcal{D}}
 
 \newcommand{\s}{\mathcal{S}}
 \newcommand{\A}{\mathcal{A}}

\newcommand{\K}{\mathcal{K}}

\begin{document}
\title[Nonlinear critical problem on compact manifolds]
{Existence and multiplicity results for nonlinear critical Neumann
problem on compact Riemannian manifolds.}
\author{ Youssef MALIKI }
\address{Department of Mathematics, Univesrity Aboubekr Belka\"{i}d of Tlemcen, ALgeria}
\email{malyouc@yahoo.fr} .
\begin{abstract}
In this work, we study on compact Riemannian manifolds with
boundary, the problems of existence and multiplicity of solutions to
a Neumann problem involving the $p-$Laplacian operator and critical
Sobolev exponents.
\end{abstract}
\maketitle
\section{Introduction}
Let $(\overline M,g)$ be an $n(n\geq3)-$dimensional Riemannian
manifolds with interior $M$ and boundary $\partial M$ which is an
$(n-1)-$dimensional Riemannian manifold with induced metric $g$.\\
In this paper, we are interested in the problem of finding solutions
on $\overline M$ of the nonlinear Neumann boundary value problem
\begin{equation}\label{eq1.1}
    \left\{
       \begin{array}{ll}
         \Delta_pu+a(x)|u|^{p-2}u=f(x)|u|^{p^*-2}u& \hbox{ in $M$;}\\
         |\nabla_g u|^{p-2}\partial_{\nu_g}u+k(x)|u|^{p-2}u=K(x)|u|^{p^{**}-2}u&
 \hbox{on $\partial M$.}
       \end{array}
     \right.
\end{equation}
where $ \partial_{\nu_g}$ is the outer unit normal derivative,
$p\in(1,n),\Delta_p=-div(|\nabla_gu|^{p-2}\nabla_gu)$ is the
$p-$Laplacian
operator, $p^*=\frac{np}{n-p}$ and $p^{**}=\frac{(n-1)p}{n-p}$.\\
The case $p=2$  corresponds to the famous problem of prescribing
scalar and mean curvatures which has been studied by several
authors, we cite for instance \cite[22]{Escobar1}, \cite{cherrier},
\cite{Ambrosetti-Li-MAlchiodi},\cite{Ambrosetti-Li-MAlchiodi1},
\cite{Djadli-Malchiodi-Ahmedou}. \\
For equation \eqref{eq1.1} without boundary condition, existence and
multiplicity results on compact and complete manifolds are obtained
by some authors \cite{aubin-cotsiolis}, \cite[9,10]{BenaliliMaliki1}
and
\cite{Druet}.\\
Problems of type of \eqref{eq1.1} are studied recently by a number
of authors. For example, in \cite{cotsiolis-labropoulos}, the
authors proved an existence results on solid torus in $\mathbb{R}^3$
under the conditions that either the function $f$ is positive and
$K$ is arbitrary or $f$ is nonnegative and  $K$ is positive jointly
with some condition on the function $a$. Successively, in
\cite{cotsiolis-labropoulos1} they obtained similar result on
general Riemannian manifolds under the same conditions.\\
Independently, in \cite{Yasov-Runst}, the authors dealt with problem
\eqref{eq1.1} in the case where the function $f$ changes sing and
the function $K$ is non positive. They obtained an existence result
for the critical case and multiplicity result for the subcritical
case. They used a fibering method introduced by Pohozaev\cite{phozaev}.\\
In the present work, we are interested  in the case where both
functions $f$ and $K$ change sign on $\overline M$. Needless to say,
this situation is more complicated and shows serious difficulties
especially when we want to define a suitable constraint set.\\
In our study, we will adapt to our case some variational techniques
introduced in \cite{rauzy} where the author dealt  with a prescribed
scalar curvature equation on compact manifolds without boundary.
These techniques rely on considering minimization problem on
suitable sets from which one can construct a sequence of curves
satisfying a certain geometry that allows to get multiplicity of
solutions in the subcritical case. By imposing additional
conditions, we improve our results by proving
multiplicity of solutions for problem \eqref{eq1.1} which includes the geometric case.\\
Note that these techniques have been adapted for equation
\eqref{eq1.1} without boundary condition in \cite{BenaliliMaliki3}
and for a $Q-$curvature equation in
\cite{Benalili}.\\\\
 {\bf Acknowledgments.} I would like to express deep gratitude to
Prof. Boumedien Abdellaoui for his valuable discussion and step by
step revision of this manuscript.
\section{Notations}
In the whole of the paper we denote by,
\begin{enumerate}
  \item $a$ and $k$ two negative constants,
\item $f$ and $K$ two changing functions respectively on $M$ and $\partial
M$,
\item $q$ and $r$ two constants such that $p<q\leq p^*, p<r\leq
p^{**}$ and $r<q$,
  \item $\|u\|_{p,M},\|u\|_{p,\partial M}$ the $L_p$ norms
respectively on $M$ and $\partial M$,
\item $H_1^p(M)$ the Sobolev space of the functions in $L_p$ with gradient in $L_p$,
\item $K_1, K_2$ the best constants defined in
the Sobolev and trace Sobolev inequalities which are the best
constants such that there exist positive constants $A$ and $B$ such
that
\begin{eqnarray*}
\|u\|_{q,M}^p&\leq&K_1\|\nabla_gu\|_{p,M}^p+A\|u\|_{p,M}^p,
p<q\leq p^* \\
\|u\|_{r,\partial M}^p &\leq&K_2
\|\nabla_gu\|_{p,M}^p+B\|u\|_{p,\partial M}^p, p< r\leq p^{**}
\end{eqnarray*}
\item $h^-=\min(h,0),h^+=\max(h,0)$, respectively the negative and
positive part of a function $h$.
\end{enumerate}
By a solution of problem \eqref{eq1.1}, we mean a function $u\in
H^p_1(M)$ such that for every $v\in \C^\infty(\overline M)$ we have
\begin{eqnarray*}
  &\int_M|\nabla_gu|^{p-2}g(\nabla_gu,\nabla_gv)dv_g+
\int_Ma|u|^{p-2}uvdv_g+\int_{\partial M}k|u|^{p-2}uvd\sigma_g& \\
  &=\int_Mf|u|^{p^*-2}uvdv_g+\int_{\partial M}K|u|^{p^{**}-2}uv
d\sigma_g.&
 \end{eqnarray*}
By regularity results \cite{Lieberman},we get that
$u\in\C^{1,\beta}(\overline M)$, for some $\beta\in(0,1)$.
\section{Statement of the results}
Our purpose is to prove existence and multiplicity of solutions to
problem \eqref{eq1.1}, it is to seek critical points $u\in
H^p_1(M)$ of the energy functional
\begin{eqnarray*}
&E(u)=\int_M|\nabla_gu|^{p}dv_g+a\int_M|u|^pdv_g+k\int_{\partial
M}|u|^pd\sigma_g&\\&-\frac{n-p}{n}\int_Mf|u|^{p^*}dv_g-\frac{n-p}{n-1}\int_{\partial
M}K|u|^{p^{**}}d\sigma_g.&
\end{eqnarray*}
 Define the quantity
\begin{equation*}
\lambda_{f,K}=\inf_\A\frac{\|\nabla_gu\|_{p,M}^p}
{|a|\|u\|_{p,M}^p+|k|\|u\|_{p,\partial M}^p}
\end{equation*}
where
\begin{equation*}
\A=\{u\in H_1^p(M),u\geq0: \int_M|f^-|udv_g+\int_{\partial
M}|K^-|ud\sigma_g=0\}
\end{equation*}
In this paper, we prove the following theorems
\begin{thm}\label{thm3.1}
Let $(M,g)$ be a compact Riemannian manifold with smooth boundary
$\partial M$.\\
There exist two positive constants $N$ and $H$ such that if the
functions $f$ and $K$ satisfies the following conditions
\begin{enumerate}
\item $\lambda_{f,K}>1$,
\item $\frac{\sup_Mf}{\int_{M}|f^-|dv_g}<N, \sup_Mf>0$
\item $\frac{\sup_{\partial M}K}{\int_{M}|K^-|d\sigma_g}<H, \sup_{\partial }K>0$
\end{enumerate}
then problem \eqref{eq1.1} admits a solution.
\end{thm}
The following theorem concerns the multiplicity of the problem
\begin{equation}\label{eq3.1}
    \left\{
       \begin{array}{ll}
         \Delta_pu+a(x)|u|^{p-2}u=f(x)|u|^{q-2}u & \hbox{ in $M$;}\\
         |\nabla_g u|^{p-2}\partial_{\nu_g}u+k(x)|u|^{p-2}u=K(x)|u|^{r-2}u &
 \hbox{on $\partial M$.}
       \end{array}
     \right.
\end{equation}
where $p<q<p^*$ and $p<r<p^{**}$.
\begin{thm}\label{thm3.2}
Let $(M,g)$ be a compact Riemannian manifold with smooth boundary
$\partial M$.\\ There exist two positive constants $N$ and $H$ such
that if the functions $f$ and $K$ satisfy the following conditions
\begin{enumerate}
\item $\lambda_{f,K}>1$,
\item $\frac{\sup_Mf}{\int_{M}|f^-|dv_g}<N,\sup_Mf>0, $
\item $\frac{\sup_{\partial
M}K}{\int_{M}|K^-|d\sigma_g}<H, \sup_{\partial M}K>0$,
\end{enumerate}
then problem \eqref{eq3.1} admits at least two distinct solutions.
\end{thm}
Under further conditions we prove the multiplicity of problem
\eqref{eq1.1}. Let $u$ be the solution given by theorem \ref{thm3.1}
and $N,H$ be the constants given in theorem \ref{thm3.2}.\\
Denote $I$ the functional
\begin{equation*}
    I(u)=\int_M|\nabla_gu|^{p}dv_g+a\int_M|u|^pdv_g+k\int_{\partial
M}|u|^pd\sigma_g.
\end{equation*} We prove
the following theorem
\begin{thm}\label{thm3.3}
Let $(M,g)$ be a compact Riemannian manifold with smooth boundary
$\partial M$.\\ Suppose that the functions $f$ and $K$ satisfies the
following conditions
\begin{enumerate}
\item $\lambda_{f,K}>1$,
\item $\frac{\sup_Mf}{\int_{M}|f^-|dv_g}<N,\sup_Mf>0$ ,
\item $\frac{\sup_{\partial
 M}K}{\int_{M}|K^-|d\sigma_g}<H, \sup_{\partial M}K>0$,
\item $\frac{p(n-1)}{n(p-1)}(\sup_Mf)^{-\frac{n-p}{p}}(\sup_{\partial
M}K)^\frac{n-p}{p-1}K_1^{-n}K_2^{\frac{p(n-1)}{p-1}}\le1 $
\end{enumerate}
and that there exists  a positive function $\Phi\in H^p_1(M)$such
that $I(\Phi)>0, \int_Mf\Phi^{p^*}dv_g>0, \int_{\partial
M}K\Phi^{p^{**}}d\sigma_g>0$ and
\begin{eqnarray*}
  0<\underset{\lambda\in[0,1]}\sup E(\lambda\Phi)<
E(u)+\frac{p}{n}\left[\sup_Mf\right]^{1-\frac{n}{p}}K_1^{-\frac{n}{p}}
\end{eqnarray*}
then, problem \eqref{eq1.1} admits at least two distinct solutions.
\end{thm}
\section{Necessary condition}
Recall that we have defined  $\lambda_{f,K}$ as
\begin{equation*}
    \lambda_{f,K}=\inf_\A\frac{\|\nabla_gu\|_{p,M}^p} {|a|\|u\|_{p,M}^p+|k|\|u\|_{p,\partial M}^p}
\end{equation*}
where
\begin{equation*}
    \A=\{u\in H_1^p(M),u\geq0: \int_M|f^-|udv_g+\int_{\partial M}|K^-|ud\sigma_g=0\}
\end{equation*}
Let us prove the following  lemma
\begin{lem}\label{lem3.1}
If the problem \ref{eq1.1} admits a solution, then
$\lambda_{f,K}\geq1$
\end{lem}
\begin{proof} First, we show that $\lambda_{f,K}$ is attained. By homogeneity, we may
take $\{u_i\}\subset H^p_1(M),u_i\geq0$ a minimizing sequence for
$\lambda_{f,K}$ such that $\|u_i\|_{p,M}^p+\|u_i\|_{p,\partial
M}^p=1$, then $u_i$ is bounded in $H^p_1(M)$. By the
Rellich-Kondrakov and Banach theorems, there exist a subsequence
$u_i$ and a function $u$ such that $u_i$ converges weakly in
$H^p_1(M)$, strongly in $L_s(M)$ and $L_t(\partial M),s<p^*,
t<p^{**}$, almost everywhere in $M$ and in the sens of trace on
$\partial M$. Then we get $\|u\|_{p,M}^p+\|u\|_{p,\partial M}^p=1$
and $u\in \A$. Moreover, the weak convergence gives that
\begin{equation*}
    \|\nabla_gu\|_{p,M}^p\leq\lim\inf\|\nabla_gu_i\|_{p,M}^p,
\end{equation*}
thus
\begin{equation*}
\lambda_{f,K}=\frac{\|\nabla_gu\|_{p,M}^p}
{|a|\|u\|_{p,M}^p+|k|\|u\|_{p,\partial M}^p}.
\end{equation*}
which means that $\lambda_{f,K}$ is attained by the function $u$ and
by regularity theorems $u$ is $\C^{1,\alpha}(M)$ for certain
$\alpha\in(0,1)$.\\
Now, we state the following generalized Picone's inequality
\cite{Abdellaoui-Peral}: for two differentiable functions $u\ge0$
and $v>0$, we have
\begin{equation}\label{eqn1}
|\nabla_gu|^p\ge|\nabla_gv|^{p-2}g
\left(\nabla_gv,\nabla_g\frac{u^p}{v^{p-1}}\right).
\end{equation}
Take $u$ a minimizer $\lambda_{f,K}$ and let $v$ be a positive
solution of the problem \eqref{eq1.1}, then we have
\begin{eqnarray*}
&\int_M|\nabla_gv|^{p-2}g
\left(\nabla_gv,\nabla_g\frac{u^p}{v^{p-1}}\right)dv_g=|a|\|u\|_{p,M}+|k|\|u\|_{p,\partial
M}\\&+\int_Mf^+u^pdv_g+\int_{\partial M}K^+u^pdv\sigma_g
\end{eqnarray*}
so, by \eqref{eqn1} we get
\begin{equation*}
\lambda_{f,K}=\frac{\|\nabla_gu\|_{p,M}}{|a|\|u\|_{p,M}+|k|\|u\|_{p,\partial
M}}\geq1
\end{equation*}
\end{proof}
\section{Subcritical problem: Multiplicity result}
In this section, we prove a multiplicity result for problem
\ref{eq3.1}, we look for solutions as critical points of the
following functional
\begin{eqnarray*}
    &&E_{q,r}(u)=\int_M|\nabla_gu|^pdv_g+a\int_M|u|^pdv_g+k\int_ {\partial
    M}|u|^pd\sigma_g\\&&-\frac{p}{q}\int_Mf|u|^{q}dv_g-\frac{p}{r}\int_{\partial
    M}K|u|^{r},q\in]p,p^*],r\in]p,p^{**}].
\end{eqnarray*}
which is bounded on the set
\begin{equation*}
    \s_{\ell,q,r}=\{u\in H^p_1(M):\|u\|_{q,M}^q+\|u\|_{r,\partial
    M}^r=\ell\}.
\end{equation*}
set
\begin{equation*}
    \mu_{\ell,q,r}=\inf_{\s_{\ell,q,r}}E(u)_{q,r}
\end{equation*}
We prove the following lemmas
\begin{lem}\label{lem5.1}
$\mu_{\ell,q,r}$ is attained for $q<p^*$ and $r<p^{**}$.
\end{lem}
\begin{proof}
Let $\ell>0$ and $\{u_i\}_{i\geq1}\subset\s_{\ell,m,q,r}$ be a
minimizing sequence; that is $\|u_i\|_{q,M}^q+\|u_i\|_{r,\partial
    M}^r=\ell$ and
    $\lim_{i\rightarrow\infty}E_{q,r}(u_i)=\mu_{\ell,m,q,r}$.
For $i$ large enough we can assume
$E_{q,r}(u_i)\leq\mu_{\ell,q,r}+1$, which gives that
\begin{eqnarray*}
&\|\nabla_gu_i\|_{p,M}^p\leq\mu_{\ell,q,r}
+|a|Vol(M)^{1-\frac{p}{q}}\ell^{\frac{p}{q}}+ |k|Vol(\partial
M)^{1-\frac{p}{r}}\ell^{\frac{p}{r}}&\\&+\left(\sup_Mf+\sup_{\partial
M}K\right)\ell+1.&
\end{eqnarray*}
Thus, the sequence $\{u_i\}_{i\geq1}$ is bounded in $H_1^p(M)$.
Since the inclusions of  $H_1^p(M)$ in $L_p(M),L_p(\partial
M),L_q(M)$ and $L_r(\partial M)$ are compact for $q<p^*$ and
$r<p^{**}$, then there exist a subsequence $u_i$ and a function
$u\in H_1^p(M)$ such that $u_i$ converges  to $u$ weakly in
$H_1^p(M)$ and strongly in each of the spaces $L_p(M), L_p(\partial
M)$,$ L_q(M), L_r(\partial M)$. It converges also to $u$ almost
everywhere in $M$ and in sense of trace on $\partial M$. Thus,
$u\in\s_{\ell,m,q,r}$ and $E_{q,r}(u)\geq\mu_{\ell,m,q,r}$.
Moreover, the strong and weak convergence imply that
\begin{equation*}
    \|\nabla_gu\|_{M,p}\leq\lim\inf\|\nabla_gu_i\|_{M,p}
\end{equation*}
hence
\begin{equation*}
 E_{q,r}(u)\leq\lim_{i\rightarrow\infty}E_{q,r}(u_i)=\mu_{\ell,q,r}
\end{equation*}
which means that $\mu_{\ell,m,q,r}$ is attained.
\end{proof}
\begin{lem}
$\mu_{\ell,q,r}$ is continuous as a function in the variable $\ell$.
\end{lem}
\begin{proof}
Let $\ell\in\left]0,\infty\right[$ and
$\{\ell_n\}_{n\in\mathbb{N}}\subset\left]0,\infty\right[$ be a
sequence such that $\lim_{n\rightarrow\infty}\ell_n=\ell<\infty$. By
lemma \ref{lem5.1}, for every $n\in\mathbb{N}$ there exist
$u_n\in\s_{\ell_n,q,r}$  such that $E_{q,r}(u_n)=\mu_{\ell_n,q,r}$
and there exists a function $u\in\s_{\ell,q,r}$ such that $
E_{q,r}(u)=\mu_{\ell,q,r}$. The sequence $u_n$ is bounded in the
spaces $L_q(M), L_r(\partial M),L_p(M)$, $L_p(\partial M)$ and
satisfies
\begin{eqnarray*}
&\|\nabla_gu_n\|_{p,M}^p\leq\mu_{\ell_n,q,r}+
|a|vol(M)^{1-\frac{p}{q}}\ell_n^{\frac{p}{q}}+|k|vol(\partial
M)^{1-\frac{p}{r}}\ell_n^{\frac{p}{r}}
&\\&+\left(\sup_Mf+\sup_{\partial M}K\right)\ell_n.
\end{eqnarray*}
On the other hand, there exits
$t_n>0,\underset{n\rightarrow\infty}{\lim}t_n=1$ such that
$t_{n}u\in\s_{\ell_n,q,r}$. Then
\begin{equation*}
\mu_{\ell_n,q,r}\leq{t_n}^p\|\nabla_gu\|_{p,M}^p+\left(t_n^q|\inf_Mf|+
t_n^r|\inf_{\partial M}K|\right)\ell,
\end{equation*}
which gives that the sequence $\{u_n\}_{n\geq0}$ is bounded in
$H_1^p(M)$. Up to a subsequence $u_n$ converges to a function
$\hat{u}$ strongly in the spaces $L_q(M), L_r(\partial M),L_p(M)$,
$L_p(\partial M)$ and weakly in $H_1^p(M)$. Thus $\hat{u}\in
\s_{\ell,q,r}$ which implies that
$E_{q,r}(\hat{u})\geq E_{q,r}(u)$.  \\
Hence,
\begin{eqnarray*}
\lim_{n\rightarrow\infty}\inf E_{q,r}(u_n)-E_{q,r}(u)\geq0.
\end{eqnarray*}
On the other hand, we have
\begin{eqnarray*}
E_{q,r}(u_n) &\leq& t_n^p\left(\|\nabla_gu\|_{p,M}^p+a\|u\|_{p,M}^p+
k\|u\|_{p,\partial M}^p\right) \\
&&-t_n^q\frac{p}{q}\int_Mf|u|^qdv_g-t_n\frac{p}{r}\int_{\partial
 M}K|u|^rd\sigma_g,
 \end{eqnarray*}
thus
\begin{equation*}
    \lim_{n\rightarrow\infty}\sup E_{q,r}(u_n)-E(u)_{q,r}\leq0.
\end{equation*}
Therefore, $\mu_{\ell,m,q,r}$ is continuous.
\end{proof}
\begin{lem}
$\mu_{\ell,q,r}$ is negative for $l$ small.
\end{lem}
\begin{proof}
For $\ell>0$, let $m_\ell>0$ be the solution of the equation
$vol(M)m_\ell^q+vol(\partial M)m_\ell^r=\ell$. Take the constant
function $u=m_\ell$, then $u\in\s_{\ell,q,r}$ and
\begin{eqnarray*}
\mu_{\ell,q,r}\leq E_{q,r}(u)&=&
m_\ell^{p}\left[-|a|vol(M)-|k|vol(\partial
M)-\frac{p}{q}m_\ell^{q-p}\int_M fdv_g\right.\\&&
\left.-\frac{p}{r}m_\ell^{r-p}\int_{\partial M}Kd\sigma_g\right],
\end{eqnarray*}
If $\ell$ is small, then so is $m_\ell$ and we conclude the lemma.
\end{proof}
\begin{lem}$\mu_{\ell,q,r}$ is negative for $\ell$ big.
\end{lem}
\begin{proof}
Let $u$ and $v$  be respectively smooth functions defined on $M$ and
$\partial M$ with supports included respectively in the sets where
$f(x)>0$ and $K(x)>0$ and such that $\|u\|_{q,M}^q=1$ and
$\|v\|_{r,\partial M}^r=1$. For $x\in \overline M$, set
$h(x)=\frac{1}{2}(\ell^\frac{1}{q}u(x)+\ell^\frac{1}{r}v(x))$, then
$h\in\s_{\ell,q,r}$ and we get\\
\begin{eqnarray*}
\mu_{\ell,q,r}\leq E(h)=\frac{1}{2^p}\ell\left[
\left(\|\nabla_gu\|^p_{p,M}+a\|u\|_{p,M}^p\right)\ell^{\frac{p}{q}-1}+
k\|u\|_{p,\partial M}^p\ell^{\frac{p}{r}-1}-\right.&&\\ \left.
\frac{1}{2^{q-p}}\int_Mf|u|^qdv_g-\frac{1}{2^{r-p}}\int_{\partial
M}K|v|^rd\sigma_g\right].&&
\end{eqnarray*}
since $\int_Mf|u|^qdv_g>0$ and $\int_{\partial M}K|v|^rd\sigma_g>0$,
then for $\ell$ large enough $\mu_{\ell,q,r}<0$.
\end{proof}
Now, let us define the quantity
\begin{equation*}
\lambda_{f,K,\eta,q,r}=\inf_{\K(\eta,q,r)}
\frac{\|\nabla_gu\|^p_{p,M}}{|a|\|u\|_{p,M}^p+|k|\|u\|_{p,\partial
M}^p}
\end{equation*}
where
\begin{eqnarray*}
&&\K(\eta,q,r)=\left\{u\in H_1^p(M):\|u\|_{q,M}^q+\|u\|_{r,\partial
M}^r=1 \text{ and }\right.\\&& \int_M|f^-||u|^qdv_g+\int_{\partial
M}|K^-||u|^rd\sigma_g \left.\leq\eta\left(\int_M|f^-|dv_g
+\int_{\partial M}|K^-|d\sigma_g\right)\right\}.
\end{eqnarray*}
\begin{lem}\label{lem5.5}
 For $(q,r)\in]p,p^*[\times]p,p^{**}[$,
$\lambda_{f,K,\eta,q,r}$ converges to $\lambda_{f,K}$ when $\eta$
goes to zero.
\end{lem}
\begin{proof}
First let us prove that $\lambda_{f,K,\eta,q,r}$ is attained. Let
$u_{i,\eta}\in\K(\eta,q,r)$ be a minimizing sequence for
$\lambda_{f,K,\eta,q,r}$, then for $i$ large enough we can have
\begin{equation*}
\|\nabla_gu_{i,\eta}\|^p_{p,M}
\leq\lambda_{f,K,\eta,q,r}(|a|{\|u_{i,\eta}\|_{p,M}^p+|k|\|u_{i,\eta}\|_{p,\partial
M}^p})+1.
\end{equation*}
Hence, the sequence $u_{i,\eta}\in\K(\eta,q,r)$ is bounded in
$H_1^p(M)$ and up to a subsequence $u_{i,\eta}$ converges to a
function $u_{\eta}$ weakly in $H_1^p(M)$ and strongly  in each of
the spaces $L_p(M),L_p(\partial M),L_q(M),L_r(\partial M)$. Thus,
$u_{\eta}\in\K(\eta,q,r)$. Moreover,
\begin{eqnarray*}
\|\nabla_gu_{\eta}\|_{M,p}\leq \lim_{i\rightarrow\infty}\inf
\|\nabla_gu_{i,\eta}\|_{M,p}.
\end{eqnarray*}
Hence
\begin{equation*}
 \frac{ \|\nabla_gu_{\eta}\|_{M,p}}
{{|a|\|u_{\eta}\|_{p,M}^p+|k|\|u_{\eta}\|_{p,\partial
M}^p}}=\lambda_{f,K,\eta,q,r}.
\end{equation*}
Now, consider $ u_{\eta}$ a sequence of $\eta$. First we observe
that if $u\in\A$ then there exists $\beta>0$ such that $\beta
u\in\K(\eta,q,r)$. By homogeneity, we get then that the sequence
$\lambda_{f,K,\eta,q,r}$ is bounded by $\lambda_{f,K}$. Thus, the
sequence $ u_{\eta}$ is bounded in $H_1^p(M)$ and there exists a
subsequence that converges, when $\eta$ goes to zero, to a function
$u$ weakly in $H_1^p(M)$ and strongly in $L_p(M),L_p(\partial
M),L_q(M)$ and
$L_r(\partial M)$. Hence, $u\in \A$.\\
On the other hand
\begin{equation*}
\|\nabla_gu\|^p_{p,M}\leq\lim_{\eta\rightarrow0}\inf
\|\nabla_gu_{\eta}\|^p_{p,M}.
\end{equation*}
Thus
\begin{equation*}
\lambda_{f,K}\leq\frac{\|\nabla_gu\|^p_{p,M}}{|a|\|u\|_{p,M}^p+|k|\|u\|_{p,\partial
M}^p}\leq\lim_{\eta\rightarrow0}\inf\lambda_{f,K,\eta,q,r}.
\end{equation*}
\end{proof}
Now, coming back to the function
$\ell\in(0,\infty)\longrightarrow\mu_{\ell,q,r}$ to prove the
following lemma
\begin{lem}\label{lem5.6} Suppose that $\lambda_{f,K}>1$, then there exist
two positive constants $N_{q}$ and  $H_{r}$ such that if
$\frac{\sup_Mf}{\int_M|f^-|dv_g}\leq N_{q}$ and
$\frac{\sup_{\partial M}K}{\int_{\partial M}|K^-|d\sigma_g}\leq
H_{r}$ then there exists an  interval $\left[\ell_1,\ell_2\right]$
such that $\mu_{\ell,q,r}$ is positive for every $\ell\in
\left[\ell_1,\ell_2\right]$.
\end{lem}
\begin{proof}First, if $\lambda_{f,K}>1$, it follows from
lemma \ref{lem5.5} that for $\eta$
small enough $\lambda_{f,K,\eta,q,r}-1=\frac{\epsilon}{|a|}>0$.\\
Let $u\in H^p_1(M)$ be such that $\|u\|_{q,M}^q+\|u\|_{r,\partial
M}^r=\ell$  with \begin{equation}\label{eq5.1} \ell\geq z=\max
\left(\left(\frac{2vol(M)^{1-\frac{p}{q}}|a|}
{\eta\int_M|f^-|dv_g}\right)^{\frac{q}{q-p}},\left(
\frac{2vol(\partial M)^{1-\frac{p}{r}}|k|}{\eta\int_{\partial
M}|K^-|d\sigma_g}\right)^{\frac{r}{r-p}}\right),
\end{equation}
take $\eta>0$ small enough so that $\ell>1$.\\
Denote by $G_{q,r}$ the functional
\begin{eqnarray*}
   &G_{q,r}(u)=\|\nabla_gu\|^p_{p,M}+a\|u\|^p_{p,M}+k\|u\|^p_{p,\partial
    M}\\&+\frac{p}{q}\int_M|f^-||u|^qdv_g+
\frac{p}{r}\int_{\partial M}|K^-||u|^rd\sigma_g.
\end{eqnarray*}
We distinguish to cases:\\
either
\begin{eqnarray*}
\frac{p}{q}\int_M|f^-||u|^qdv_g+\frac{p}{r}\int_{\partial
M}|K^-||u|^rd\sigma_g\geq\eta\ell(\int_M|f^-|dv_g +\int_{\partial
M}|K^-|d\sigma_g)
\end{eqnarray*}
and so
\begin{eqnarray*}
G_{q,r}(u)&\geq&
vol(M)^{1-\frac{p}{q}}a\ell^{\frac{p}{q}}+vol(\partial
M)^{1-\frac{p}{r}}k\ell^{\frac{p}{r}}+\ell(\eta\int_M|f^-|dv_g
+\eta\int_{\partial
M}|K^-|d\sigma_g)\\&\geq&vol(M)^{1-\frac{p}{q}}|a|\ell^{\frac{p}{q}}
\left[\frac{\eta\int_M|f^-|dv_g}{vol(M)^{1-\frac{p}{q}}|a|}\ell^{1-\frac{p}{q}}-1\right]+
vol(\partial M)^{1-\frac{p}{r}}|k|\ell^{\frac{p}{r}}
\\&& \left[\frac{\eta\int_{\partial M}|K^-|d\sigma_g}{vol(\partial
M)^{1-\frac{p}{r}}|k|}\ell^{1-\frac{p}{r}}-1\right]
\\&\geq&vol(M)^{1-\frac{p}{q}}|a|\ell^{\frac{p}{q}}+
vol(\partial M)^{1-\frac{p}{r}}|k|\ell^{\frac{p}{r}}.
\end{eqnarray*}
Or
\begin{eqnarray*} \frac{p}{q}\int_M|f^-||u|^qdv_g+\frac{p}{r}\int_{\partial
M}|K^-||u|^rd\sigma_g\leq\ell\eta(\int_M|f^-|dv_g +\int_{\partial
M}|K^-|d\sigma_g).
\end{eqnarray*}
In this case , let $\delta>0$ be a solution of the equation
$\delta^q \|u\|_{q,M}^q+\delta^r\| u\|_{r,M}^r=1$, it can be easily
seen that $\ell^{-\frac{1}{r}}<\delta<\ell^{-\frac{1}{q}}$. This
implies that $\delta u\in\K(\eta,q,r)$. In particular,
\begin{equation*}
\|\nabla_gu\|^p_{p,M}\geq\lambda_{f,K,\eta,q,r}\left[|a|\|u\|^p_{p,M}
+|k|\|u\|^p_{p,\partial
    M}\right],
\end{equation*}
so
\begin{eqnarray*}
G_{q,r}(u) &\geq& (\lambda_{f,K,\eta,q,r}-1)\left[|a|\|u\|^p_{p,M}
+|k|\|u\|^p_{p,\partial
    M}\right],
\end{eqnarray*}
write
$\min(|a|,|k|)(\lambda_{f,K,\eta,q,r}-1)=\alpha_\eta+\beta_\eta$
such that $\frac{\alpha_\eta(|a|+|k|)}{\beta}=\frac{A+B}{K_1+K_2}$,
where $A,B,K_1,K_2$ are the constants appearing in the Sobolev and
trace Sobolev inequalities. Then
\begin{eqnarray*}
G_{q,r}(u) &\geq\alpha_\eta(\|u\|^p_{p,M}+\|u\|^p_{p,\partial
M})+\frac{\beta}{|a|+|k|}\left[-
G_{q,r}(u)+\|\nabla_gu\|_{p,M}^p\right.\\&
\left.+\frac{p}{q}\int_M|f^-||u|^qdv_g +\frac{p}{r}\int_{\partial
M}|K^-||u|^rd\sigma_g\right]
\end{eqnarray*}
that is
\begin{eqnarray*}
\left(1+\frac{\beta_\eta}{|a|+|k|}\right)G_{q,r}(u)
&\geq\frac{\beta_\eta}{|a|+|k|}\left[
\|\nabla_gu\|_{p,M}^p+\frac{\alpha_\eta(|a|+|k|)}{\beta_\eta}(\|u\|^p_{p,M}+
\|u\|^p_{p,\partial M})\right].
\end{eqnarray*}
Then, by the Sobolev and trace sobolev inequalities we get
\begin{eqnarray*}
\left(1+\frac{\beta_\eta}{|a|+|k|}\right)G_{q,r}(u)
&\geq&\frac{\beta_\eta}{|a|+|k|}\left[
\|\nabla_gu\|_{p,M}^p+\frac{A+B}{K_1+K_2}(\|u\|^p_{p,M}+\|u\|^p_{p,\partial
M})\right]\\&\geq&\frac{\beta_\eta}{(|a|+|k|)(K_1+K_2)}(\ell^{\frac{p}{q}}+\ell^{\frac{p}{r}}).
\end{eqnarray*}
Thus
\begin{eqnarray*}
G(u)\geq\frac{\beta_\eta}
{(K_1+K_2)(|a|+|k|+\beta_\eta)}(\ell^{\frac{p}{q}}+\ell^{\frac{p}{r}}).
\end{eqnarray*}
On the other hand, we have
\begin{eqnarray*}
    E_{q,r}(u)&=&G_{q,r}(u)-\int_Mf^+|u|^qdv_g-\int_{\partial
    M}K^+|u|^rd\sigma\\&\geq&G_{q,r}(u)-(\sup_Mf+\sup_{\partial
    M}K)\ell\\&\geq& t\ell^{\frac{p}{q}}+s\ell^{\frac{p}{r}}-(\sup_Mf+\sup_{\partial
    M}K)\ell,
\end{eqnarray*}
where
\begin{equation*}
t=\min\left(\frac{\beta_\eta} {(K_1+K_2)(|a|+|k|+\beta_\eta)},
|a|vol(M)^{1-\frac{p}{q}}\right),
\end{equation*}
and
\begin{equation*}
s=\min\left(\frac{\beta_\eta} {(K_1+K_2)(|a|+|k|+\beta_\eta)},
|k|vol(\partial M)^{1-\frac{p}{r}}\right).
\end{equation*}
Now, the function $g:\ell\longrightarrow
\frac{t}{2}\ell^{\frac{p}{q}}+\frac{s}{2}\ell^{\frac{p}{r}}-(\sup_Mf+\sup_{\partial
M}K)\ell$ attains a positive maximum in the interval $]0,\ell_o[$
where $g(\ell_o)=0$. Therefore, the hypotheses on the functions $f$
and $K$ assure that $g(z)>0$, where $z$ is defined by \eqref{eq5.1}.
In fact,
\begin{equation*}
    g(z)\geq z\left[\frac{t\eta\int_M|f^-|dv_g}
{4vol(M)^{1-\frac{p}{q}}|a|}-\sup_Mf+\frac{s\eta\int_{\partial
M}|K^-|dv_g} {4vol(\partial M)^{1-\frac{p}{r}}|k|}-\sup_{\partial
M}K\right]
\end{equation*}
so if we require that
\begin{equation*}
    \frac{\sup_Mf}{\int_M|f^-|dv_g}\leq N_q=\frac{t\eta}{8vol(M)^{1-\frac{p}{q}}|a|}
\end{equation*}
and
\begin{equation*}
\frac{\sup_{\partial M}K}{\int_{\partial M}|K^-|dv_g}\leq
H_r=\frac{s\eta}{8vol(\partial M)^{1-\frac{p}{r}}|k|}
\end{equation*}
we get $g(z)>0$ and then $z<\ell_o$. Thus if
$\ell\in\left[z,\ell_o\right]$, then
\begin{equation*}
    E_{q,r}(u)>\frac{1}{2}(t\ell^{\frac{p}{q}}+s\ell^{\frac{p}{r}})>0,
\end{equation*}
and so $\mu_{\ell,q,r}>0$  for all
$\ell\in\left[\ell_1,\ell_2\right]=\left[z,\ell_o\right]$.
\end{proof}
Now, by mean of the mountain pass theorem, we show that the
existence of the interval $\left[\ell_1,\ell_2\right]$ leads to the
existence of a second critical point of the functional $E_{q,r}$.
First, we prove the following lemma
\begin{lem}
The Plais-Smale condition is satisfied for the functional $E_{q,r}$
, $q<p^*$ and $ r<p^{**}$.
\end{lem}
\begin{proof} First, we claim that each Plais-Smale sequence for the functional $E_{q,r}$
is bounded in $H^p_1(M)$. In fact, let $u_n\in H^p_1(M)$ be a
sequence such that\\
$E_{q,r}(u_n)\underset{n\rightarrow\infty}\rightarrow \gamma$ and $
E'_{q,r}(u_n)\underset{n\rightarrow\infty}\rightarrow0$, then we
have
\begin{eqnarray*}
 E_{q,r}(u_n)-\frac{1}{q}E'_{q,r}(u_n)u_n&=
 \left(1-\frac{p}{q}\right)\left[\|\nabla_gu_n\|_{p,M}^p+a\|u_n\|_{p,M}^p+
k\|u_n\|_{p,M}^p\right]
  &\\
  &-\left(\frac{1}{r}-\frac{p}{q}\right)\int_{\partial M}K|u_n|^rd\sigma_g&
  \end{eqnarray*}
  and
\begin{eqnarray*}
E_{q,r}(u_n)-\frac{1}{p}E'_{q,r}(u_n)u_n=
\left(1-\frac{p}{q}\right)\int_Mf|u_n|^qdv_g+
\left(1-\frac{p}{r}\right)\int_{\partial M}K|u_n|^rd\sigma_g
\end{eqnarray*}
so, for every $\varepsilon>0$, there exists $n_o$ such that for all
$n\geq n_o$, we have
\begin{eqnarray*}
  \left|\left(1-\frac{p}{q}\right)\int_Mf|u_n|^qdv_g+
\left(1-\frac{p}{q}\right)\int_{\partial
M}K|u_n|^rd\sigma_g-\gamma\right| &\leq&
\varepsilon+o(\|u_n\|_{H^p_1(M)})
\end{eqnarray*}
and
\begin{eqnarray*}
&\left|\left(1-\frac{p}{q}\right)
\left[\|\nabla_gu_n\|_{p,M}^p+a\|u_n\|_{p,M}^p+k\|u_n\|_{p,M}^p\right]
-\left(\frac{1}{r}-\frac{p}{r}\right)\int_{\partial
M}K|u_n|^rd\sigma_g-\gamma\right|&\\&\leq\varepsilon+o(\|u_n\|_{H^p_1(M)})&
    \end{eqnarray*}
Let $\ell>0$ be such that $\mu_{\ell,q,r}>0$ and put $v_n=\beta_{n}
u_n$, where $\beta_n>0$ is such that
$\beta_{n}^q\|u\|_{q,M}^q+\beta_{n}^r\|u\|_{r,M}^r=\ell$. We observe
that
$v_n$ is bounded in $L_q(M)$ and $L_{r}(\partial M)$.\\
Now, we have
\begin{eqnarray}\label{eq6.1}
 \nonumber
&\left|\frac{1}{\beta_n^q}\left(1-\frac{p}{q}\right)\int_Mf|v_n|^qdv_g+
\frac{1}{\beta_n^r}\left(1-\frac{p}{r}\right)\int_{\partial
M}K|v_n|^rd\sigma_g-\gamma\right|\\&\leq \varepsilon
+o\|v_n\|_{H^p_1(M)}&
\end{eqnarray}
and
\begin{eqnarray}\label{eq6.2}
  \nonumber
&\left|\frac{1}{\beta_n^p}\left(1-\frac{p}{q}\right)
\left[\|\nabla_gv_n\|_{p,M}^p+a\|v_n\|_{p,M}^p+k
\|v_n\|_{p,M}^p\right]\right.\\&\left.
-\frac{1}{\beta_n^r}\left(\frac{1}{r}-\frac{p}{q}\right)
\int_{\partial M}K|v_n|^rd\sigma_g-\gamma\right|\leq\varepsilon+
o\|v_n\|_{H^p_1(M)}.&
\end{eqnarray}
Since the sequence $v_n $ is bounded in $L_q(M)$ and $L_r(\partial
M)$, then by \eqref{eq6.2} it is bounded in $H_1^p(M)$. Moreover, we
affirm that the sequence $u_n$ is bounded in $L_q(M)$ and
$L_r(\partial M)$. In fact, if the sequence $u_n$ goes to infinity
in $L_q(M)$ or $L_r(\partial M) $ then the sequence $\beta_n$ goes
to zero when $n$ goes to infinity. This implies by mean of
inequalities \eqref{eq6.1} and \eqref{eq6.2} that
$E_{q,r}(v_n)\rightarrow0$ as $n\rightarrow\infty$.\\
Since $v_n\in\s_{\ell,m,q,r}$, then
$E_{q,r}(v_n)\geq\mu_{\ell,q,r}>0$, this is a patent contradiction.
 Thus, $u_n$ is bounded in $L_q(M)$ and $L_r(\partial M)$
 and since $E'_{q,r}(u_n)u_n\longrightarrow0$ then, it is bounded in $H^p_1(M)$.\\
Thus, up to a subsequence $u_n$ converges  to a function $u$ weakly
 in $H^p_1(M)$ and  strongly in $L_q(M),L_p(M),L_r(\partial M)$ and
$L_p(\partial M)$. Then, by Bresis-Lieb lemma, we obtain
\begin{eqnarray*}
\|\nabla_g(u_n-u)\|_{p,M}&\leq&|a|\|u_n-u\|_{p,M}^p+|k|\|u_n-u\|_{p,\partial
M}^p+\sup_Mf\|u_n-u\|_{q,M}^q\\ &&+\sup_{\partial
M}K\|u_n-u\|_{r,\partial M}^r+o(1)\\&\leq& o(1).
\end{eqnarray*}
which means that the subsequence $u_n$ converges strongly to $u$ in
$H^p_1(M)$.
\end{proof}
Now, we prove theorem \ref{thm3.2}
\begin{proof}[\textbf{Proof of theorem \ref{thm3.2}}]
\begin{quote}
\end{quote}
\begin{flushleft}
$\bullet$ \textbf{Existence of first solution}
\end{flushleft}Let $\ell_1>0$ be such that $\mu_{\ell_1,q,r}=0$ and the curve
$\ell\rightarrow\mu_{\ell,q,r}$ is negative for
$\ell\in]0,\ell_1[ $ \\
Set
\begin{equation*}
\mu_{\ell_{q,r},q,r}=\inf_{\D_{\ell,q,r}}E_{q,r}
\end{equation*}
where
\begin{equation*}
\D_{\ell,q,r}=\{u\in H^p_1(M): u>0, \|u\|_{q,M}^q+\|u\|_{r,\partial
M}^r\leq\ell,\ell<\ell_1\}
\end{equation*}
Take $\ell$ as small as $\mu_{\ell,q,r}<0$, then there exists $u\in
\D_{\ell,q,r}$ such that $\|u\|_{q,M}^q+\|u\|_{r,\partial M}^r=\ell$
and $ E_{q,r}(u)=\mu_{\ell,q,r}$ in such way that,
$u\in\D_{\ell,q,r}$ and
\begin{equation}\label{eq5.2}
\mu_{\ell_{q,r},q,r}\leq E_{q,r}(u)=\mu_{\ell,q,r}<0.
\end{equation}
By using the Ekeland Variational Principle, in the set
$\D_{\ell,q,r}$ we can find a sequence $\{u_{q,r,n}\}_n$ such that
$E_{q,r}(u_{q,r,n})\underset{
n\rightarrow\infty}\rightarrow\mu_{\ell_{q,r},q,r}$ and
$E'_{q,r}(u_{q,r,n})\underset{ n\rightarrow\infty}\rightarrow0$.
Obviously, the sequence $u_{q,r,n}$ is bounded in $H^p_1(M)$, then
up to a subsequence, $u_{q,r,n}$ converges to a function $u_{q,r}$
weakly in $H^p_1(M)$ and strongly in $L_s(M), L_t(\partial
M)(s<p^*,t<p^{**})$ , almost everywhere on $M$ and in the sense of
trace on $\partial M$. Thus, $u_{q,r}\in\D_{\ell,q,r}$ and
\begin{equation*}
 E_{q,r}(u_{q,r})\le\lim_{n\rightarrow\infty}\inf
E_{q,r}(u_{q,r,n})=\mu_{\ell_{q,r},q,r}<0.
\end{equation*}
Furthermore, by the weak convergence in $H^p_1(M),
L_{\frac{q}{q-1}}(M)$ and $L_{\frac{r}{r-1}}(\partial M)$, it
follows that for every $v\in H^p_1(M)$
\begin{equation*}
\langle E'(u_{q,r}),v\rangle=\lim_{n\rightarrow\infty}\langle
E'(u_{q,r,n}),v\rangle=0.
\end{equation*}
Hence, $u_{q,r}$ is critical point of $E_{q,r}$ with
$E_{q,r}(u_{q,r})<0$.
\begin{flushleft}
    $\bullet$ \textbf{Existence of second solution}
\end{flushleft}
Now, we prove that there exists a
second solution $v_{q,r}$  with $E_{q,r}(v_{q,r})>0$.\\
 Let
$\ell_1$ and $\ell_2$ be such that
\begin{eqnarray*}
\mu_{\ell_1,q,r}=E_{q,r}(u_{\ell_1 })&=&0\\
\mu_{\ell_2,q,r}=E_{q,r}(u_{\ell_2} )&=&0
\end{eqnarray*}
and consider
\begin{equation*}
    \mu_{q,r}=\inf_{g\in\Gamma}\max_{s\in[0,1]}E_{q,r}(g(s))
\end{equation*}
where
\begin{equation*}
\Gamma=\{g\in\C([0,1],H^p_1(M)):g(0)=u_{\ell_1 },g(1)=u_{\ell_2}\}
\end{equation*}
We claim that $\mu_{q,r}$ is critical value of the functional
$E_{q,r}$. In fact, if it is not, then there exists $\varepsilon>0$
small such that $E_{q,r}$ does not possess any critical value in the
interval $[\mu_{q,r}-\varepsilon,\mu_{q,r}+\varepsilon]$. Thus, by
the deformations Lemma we can find a function
$\phi_t:H^p_1(M)\longrightarrow H^p_1(M), t\in[0,1]$, continuous in
$t$ such that:\begin{enumerate}
                \item $\phi_0(u)=u,\forall u\in H^p_1(M)$
                \item $\phi_t(u)=u,\forall u\in H^p_1(M)$ such that
$E_{q,r}(u)\not\in[\mu_{q,r}-\varepsilon,\mu_{q,r}+\varepsilon]$
                \item $\forall u\in H^p_1(M),$ s.t $
E_{q,r}(u)\leq\mu_{q,r}+\varepsilon$, then $
 E_{q,r}(\phi_1(u))\leq\mu_{q,r}-\varepsilon$.
              \end{enumerate}
Now, let $g\in \Gamma$ be such that
$\max_{s\in[0,1]}E_{q,r}(g(s))\leq\mu_{q,r}+\varepsilon$. By
definition of $u_{\ell_1},u_{\ell_2}$ and property (2) of the
function $\phi_t$, we have
\begin{eqnarray*}
\phi_t(u_{\ell_1}) &=& u_{\ell_1}\\
 \phi_t(u_{\ell_2})&=&u_{\ell_2}
\end{eqnarray*}
In particular, the curve $\phi_1(g)\in\Gamma$ which gives that
$\mu_{q,r}\leq \max_{s\in[0,1]}E_{q,r}(\phi_1(g(s)))$, but by
property (3), we have that
$\max_{s\in[0,1]}E_{q,r}(\phi_1(g(s)))\leq\mu_{q,r}-\varepsilon$,
which makes  a contradiction. $\mu_{q,r}$ is then a critical level
for the functional $E_{q,r}$ and
\begin{equation*}
\mu_{q,r}>\sup_{\ell\in[\ell_1,\ell_2]}\mu_{\ell,q,r}>0
\end{equation*}
Therefore, there exists $v_{q,r}$ a solution  of \eqref{eq3.1} with
$E_{q,r}(v_{q,r})=\mu_{q,r}>0$ and theorem \ref{thm3.2} is proven.
\end{proof}
\section{critical problem: multiplicity result}
In this section, we prove existence and multiplicity of solutions of
the problem \eqref{eq1.1}. We will study the limits of the sequences
$u_{q,r}$ and $v_{q,r}$ as $(q,r)$ goes to$(p^*,p^{**})$. Here,
besides the non-compactness of the inclusions
$H^p_1(M)\hookrightarrow L_{p^*}(M)$ and $L_{p^{**}}(\partial M)$,
due to the fact that the functions $f$ and $K$ change the sign, we
face serious problems in proving distinction among the limits . The
curve $\ell\to \mu_{\ell,q,r}$ will play an important role in
overcoming these problems. \\
Let $E$ be the functional
\begin{eqnarray*}
&E=\int_M|\nabla_gu|^{p}dv_g+a\int_M|u|^pdv_g+k\int_{\partial
M}|u|^pd\sigma_g&\\&-\frac{n-p}{n}\int_Mf|u|^{p^*}dv_g-\frac{n-p}{n-1}\int_{\partial
M}K|u|^{p^{**}}d\sigma_g,&
\end{eqnarray*}
\subsection{Existence of first solution}
\begin{proof}[\textbf{Proof of theorem \ref{thm3.1}}]
Let $u_{q,r}>0$  be the sequence of critical points of the
functional $E_{q,r}$ such that $E_{q,r}(u_{q,r})<0$. This sequence
is bounded in $H^p_1(M)$. In fact, we have
\begin{eqnarray*}
 &\|u_{q,r}\|_{q,M}^q+\|u_{q,r}\|_{r,\partial M}^r= \ell_{q,r}<\ell_1&\\&
<\max\left(\left(\frac{2vol(M)^{1-\frac{p}{n}}|a|}
{\eta\int_M|f^-|dv_g}\right)^{\frac{n}{p}},\left(
\frac{2vol(\partial M)^{1-\frac{p}{n-1}}|k|}{\eta\int_{\partial
M}|K^-|d\sigma_g}\right)^{\frac{n-1}{p}}\right)+\varepsilon&
\end{eqnarray*}
thus $u_{q,r}$ is bounded in $L_q(M)$ and $L_r(\partial M)$.
Moreover, by  \eqref{eq5.2} we have
\begin{eqnarray*}
&\|\nabla_gu_{q,r}\|_{p,M}\leq\mu_{r,q}+|a|vol(M)^{\frac{p}{n}}\ell_1^{\frac{n-p}{n}}
+|k|vol(\partial
M)^{\frac{p-1}{n-1}}\ell_1^{\frac{n-p}{n-1}}&\\&+(\sup_{M}f+\sup_{\partial
M}K)\ell_1.&
\end{eqnarray*}
Since $\mu_{\ell,q,r}<0$, the sequence $u_{q,r}$ is bounded in
$H^p_1(M)$. Then, we can obtain a subsequence $u_{q,r}$ and a
function  $u\in H^p_1(M) $ such that
\begin{enumerate}
  \item $u_{q,r}$ converges weakly to $u$ in $H^p_1(M)$
  \item $u_{q,r}$ converges strongly to $u$ in $L_p(M), L_p(\partial M),
 L_{p^*-1}(M) $ and $L_{p^{**}-1}(\partial
M)$
  \item $u_{q,r}$ converges almost everywhere to $u$ in $M$ and in
the sens of trace on $\partial M$.
\item The sequence $\nabla_gu_{q,r}$ converges almost everywhere to $\nabla_gu$
\end{enumerate}
 Thus, $u$ is critical point of $E$. By Bresis-Lieb lemma \cite{Bresis-Lieb},
Sobolev and trace Sobolev inequalities, we get
\begin{eqnarray*}
E_{q,r}(u_{q,r})-E_{q,r}(u)&=&\|\nabla_g(u_{q,r}-u)\|^{p}_{p,M}-(1-\frac{p}{q})\int_Mf|u_{q,r}-u|^qdv_g
\\&&-(1-\frac{p}{q})\int_{\partial
M}K|u_{q,r}-u|^rd\sigma_g+o(1)\\&\geq&
\frac{1}{K_1+K_2}\left[\|u_{q,r}-u\|^{p}_{q,M}+
\|u_{q,r}-u\|^{p}_{r,M}\right]\\&&-\sup_Mf\|u_{q,r}-u\|^{q}_{q,M}-\sup_{\partial
M}K\|u_{q,r}-u\|^{r}_{r,\partial M}+o(1)
\end{eqnarray*}
then, taking into account, by Bresis-Lieb Lemma again that
$\lim_{(q,r)\rightarrow(p^*,p^{**})}
\|u_{q,r}-u\|_{q,M}^q+\|u_{q,r}-u\|^{r}_{r,\partial M}\leq\ell,
\ell>0$, we get
\begin{eqnarray*}
E(u)&\le&\lim_{(q,r)\rightarrow(p^*,p^{**})}
E_{q,r}(u_{q,r})-\frac{1}{K_1+K_2}\lim_{(q,r)\rightarrow(p^*,p^{**})}
\left[\|u_{q,r}-u\|^{p}_{q,M}\right.\\&&\left.+
\|u_{q,r}-u\|^{p}_{r,M}\right]+(\sup_Mf+\sup_{\partial M}K)\ell.
\end{eqnarray*}
Sine, a priori, we have that
$\lim_{(q,r)\rightarrow(p^*,p^{**})}\left[\|u_{q,r}-u\|^{p}_{q,M}+
\|u_{q,r}-u\|^{p}_{r,M}\right]>0 $, by taking $\ell$ small enough,
we obtain $E(u)<0$, thus $u\neq0$ and we are done.
\end{proof}
\subsection{Existence of second solution}
Now, we consider the sequence $\{v_{q,r}\}$ of solutions of the
subcritical problem \eqref{eq3.1} obtained by the mountain pass
lemma. We prove that $\{v_{q,r}\}$ will converge to non zero and
different critical point of $E$. We recall that the sequence
$\{v_{q,r}\}$ fulfills the following properties
\begin{equation} E'_{q,r}(v_{q,r})=0 \text{ and }
E_{q,r}(v_{q,r})=\mu_{q,r}>0.
\end{equation}
First, we prove the following lemma
\begin{lem}\label{lem6.1}
The sequence of functions $v_{q,r}$ is bounded in $H^p_1(M)$.
\end{lem}
\begin{proof}
Let $\ell_1>0$ and $\ell_2>0,\ell_1<\ell_2$ be two real numbers with
the associated functions $u_{\ell_1}$ and $u_{\ell_1}$ such that
\begin{eqnarray*}
\mu_{\ell_1,q,r}=E_{q,r}(u_{\ell_1 })&=&0\\
\mu_{\ell_2,q,r}=E_{q,r}(u_{\ell_2} )&=&0.
\end{eqnarray*}
For $s\in[0,1]$, let $g$ be the curve $g(s)=su_{\ell_1
}+(1-s)u_{\ell_2}$. Then by definition of $\mu_{q,r}$ we have
\begin{equation*}
    \mu_{q,r}\le\max_{s\in[0,1]}E_{q,r}(g(s)).
\end{equation*}
Since $E_{q,r}(g(1))=0$ and we can find $s\in[0,1]$ such that
$$\|g(s)\|_{q,M}^q+\|g(s)\|_{r,\partial
M}^r=\frac{\ell_1+\ell_2}{2}>\ell_1,$$  it follows that the curve
$E_{q,r}(g(s))$ attains for certain $s_o\in(0,1)$ a positive
maximum. Thus
\begin{eqnarray*}
&\mu_{q,r}\le
E_{q,r}(g(s_o))=(1-\frac{p}{q})\int_Mf|g(s_o)|^qdv_g+(1-\frac{p}{r})\int_{\partial
M}K|g(s_o)|^rd\sigma_g&\\&\le(\sup_Mf+\sup_{\partial
M}K)(\ell_1+\ell_2).
\end{eqnarray*}
which gives that the sequence $\mu_{q,r}$ is uniformly bounded in
$(q,r)$.\\
Now, it remains to show that the sequence $v_{q,r}$ is bounded in
$L_q(M)$ and $L_r(\partial M)$. We proceed as in the proof of
theorem \ref{thm3.2}. Let $\ell_o\in (\ell_1,\ell_2)$ be such that
$\mu_{\ell_o,q,r}=\sup_\ell\mu_{\ell,q,r}>0$ and consider the
sequence $u_{q,r}=\beta_{q,r}v_{q,r}$ where $\beta_{q,r}$ is such
that
$\beta_{q,r}^q\|v_{q,r}\|_{q,M}^q+\beta_{q,r}^r\|v_{q,r}\|_{r,\partial
M}^r=\ell_o$ . The sequence $u_{q,r}$, such as defined, satisfies
\begin{eqnarray}\label{eq7.2}
\mu_{q,r}=(1-\frac{p}{q})\beta_{q,r}^{-q}\int_Mf|u_{q,r}|^qdv_g+
(1-\frac{p}{r})\beta_{q,r}^{-r}\int_{\partial
M}K|u_{q,r}|^rd\sigma_g,
\end{eqnarray}
and
\begin{eqnarray}\label{eq7.3}
E_{q,r}(u_{q,r})= \beta_{q,r}^{p}\mu_{q,r}+\frac{p}{q}
(\beta_{q,r}^{p-q}-1)\int_Mf|u_{q,r}|^qdv_g\\
\nonumber+\frac{p}{r}(\beta_{q,r}^{p-r}-1)\int_{\partial
M}K|u_{q,r}|^rd\sigma_g.
\end{eqnarray}
Suppose by contradiction that the sequence $v_{q,r}$ goes to
infinity in $L_q(M)$ and $L_r(\partial M)$ as $(q,r)$ goes to
$(p^*,p^{**})$.Then, the sequence $\beta_{q,r}$ should go to zero.
Since $\mu_{q,r}$ is bounded, we get necessarily that
$\int_Mf|u_{q,r}|^qdv_g$  and $\int_{\partial
M}K|u_{q,r}|^rd\sigma_g$ go to zero  as $(q,r)$ goes to
$(p^*,p^{**})$ and $\beta_{q,r}^{-q}\int_Mf|u_{q,r}|^qdv_g$,
$\beta_{q,r}^{-r}\int_{\partial M}K|u_{q,r}|^rd\sigma_g$ are both
bounded. Thus, by \eqref{eq7.3}, we get that $E_{q,r}(u_{q,r})$ goes
to zero as as $(q,r)$ goes to $(p^*,p^{**})$.\\
On the other hand, we have $E_{q,r}>\mu_{\ell_o,q,r}$ and by lemma
\ref{lem5.6} the sequence $\mu_{\ell_o,q,r}$ does not go to zero as
$(q,r)$ goes to $(p^*,p^{**})$, this makes a contradiction.
Therefore, the sequence $v_{q,r}$ is bounded in $L_q(M$ and
$L_r(\partial M)$ and since it satisfies
\begin{eqnarray*}
\|\nabla_gv_{q,r}\|_{p,M}^p&=&
|a|\|v_{q,r}\|_{p,M}^p+|k|\|v_{q,r}\|_{p,\partial M}^p\\&+&
\int_Mf|v_{q,r}|^qdv_g+\int_{\partial M}K|u_{q,r}|^rd\sigma_g,
\end{eqnarray*}
then it is bounded in $H^p_1(M)$.
\end{proof}
Now, as the sequence $v_{q,r}$ is bounded in $H^p_1(M)$, we can
extract a subsequence that converges strongly to a function $v$ in
$L_p(M)$ and $L_p(\partial M)$ and weakly in $ H^p_1(M),
L_{\frac{p^*}{p^{*}-1}}(M)$ and
$L_{\frac{p^{**}}{p^{**}-1}}(\partial M)$. The function $v$ is then
a critical point of the functional $E$. But, this is not enough to
conclude existence of second solution because in spite of the fact
that $\underset{(q,r)\rightarrow(p^*,p^{**})}\lim
E_{q,r}(v_{q,r})-E_{q,r}(u_{q,r})>0$, we could have $v=u$ or $v=0$
regarding the lack of the strong convergence of the sequence
$v_{q,r}$ to $v$ in $L_{p^*}(M)$ and $L_{p^{**}}(\partial M)$. In
the following lemmas, we give sufficient conditions to
prevent such cases from occurrence .\\
First, note that the first solution $u$ satisfies that
$\|u\|_{p^*,M}^{p^*}+\|u\|_{p^{**},\partial
M}^{p^{**}}\leq\ell<\ell_1$. Take $\ell$ such that
\begin{equation*}
\ell<\frac{(\sup_Mf)^{-\frac{n-p}{p}}K_1^{-n}}{|\inf_Mf|+|\inf_{\partial
M}K|},
\end{equation*}
we get then
\begin{equation*}
E(u)+\frac{p}{n}\left (\sup_Mf\right)^{1-\frac{n}{p}}K_1^{-n}>0.
\end{equation*}
Let us prove the following lemma
\begin{lem}\label{lem6.2}
Suppose that the sequence $v_{q,r}$ converges strongly to the
function $u$ in $L_p(M)$ and $L_p(\partial M)$ and that the
functions $f$ and $K$ satisfy
\begin{equation*}
\frac{p(n-1)}{n(p-1)}(\sup_Mf)^{-\frac{n-p}{p}}(\sup_{\partial
M}K)^\frac{n-p}{p-1}K_1^{-n}K_2^{\frac{p(n-1)}{p-1}}\le1
\end{equation*}
If the following condition is satisfied
\begin{eqnarray}\label{eqt6.4}
\underset{(q,r)\rightarrow(p^*,p^{**})}\lim\mu_{q,r}<E(u)+\frac{p}{n}\left
(\sup_Mf\right)^{1-\frac{n}{p}}K_1^{-n}
\end{eqnarray}
then the sequence  $v_{q,r}$ converges strongly to $u$ in
$H^p_1(M)$.
\end{lem}
\begin{proof}
The sequence $v_{q,r}$ is bounded. We may assume that $v_{q,r}$
converges to $u$ weakly in $H^p_1(M)$, almost every where on $M$ and
in the sense of trace on $\partial M$.\\ In particular, we assume
that the sequences of measures
$|\nabla_gv_{q,r}|^pdv_g,|v_{q,r}|^{p^*}dv_g$ and
$|v_{q,r}|^{p^{**}}d\sigma_g$ converge weakly in the sense of
measures respectively to bounded nonnegative measures $d\mu,d\nu$ and $d\pi$.\\
Thus, by a version of a concentration-compactness theorem for
manifolds with boundary \cite{Biezuner}, there exist at most
countable index set $I$, sequence of points $\{x_i\}_{i\in
I}\subset\overline M$ and positive numbers $\{\mu_i\}_{i\in I},
\{\nu_i\}_{i\in I}, \{\pi_i\}_{i\in I}$ such that
\begin{eqnarray*}
  d\mu &\ge& |\nabla_gu|^pdv_g +\sum_{i\in I}\mu_i\delta_{x_i},\\
  d\nu &=&|u|^{p^*}dv_g+ \sum_{i\in I}\nu_i\delta_{x_i}, \text{ and } \\
  d\pi &=&|u|^{p^{**}}d\sigma_g+\sum_{i\in J}\pi_i\delta_{y_i}.
\end{eqnarray*}
Moreover, \begin{eqnarray*}
            \nu_i^{\frac{1}{p^*}}\le K_1\mu_i^{\frac{1}{p}} \text{   and }
             \pi_i^{\frac{1}{p^{**}}}\le K_2\mu_i^{\frac{1}{p}}.
          \end{eqnarray*}
Take $x_i\in\overline M$ in the support of the singular part of
$\mu,\nu$ and $\pi$ and for $\varepsilon>0$, let $\phi$ be a
$\C^{\infty}(B(x_i,\varepsilon))$ cut-off function such that
$supp\phi\subset(B(x_i,\varepsilon),\phi\equiv1$ on
$(B(x_i,\frac{\varepsilon}{2}))$ and $|\nabla_g\phi|\leq C$. Then,
we get
\begin{eqnarray*}
\int_{B(x_i,\varepsilon)}\phi
d\mu&=&\lim_{(q,r)\rightarrow(p^*,p^{**})}\int_{B(x_i,\varepsilon)}
|\nabla_gv_{q,r}^{p-2}|g
(\nabla_gv_{q,r},\nabla_g\phi)dv_g \\
&=&-a\int_{B(x_i,\varepsilon)}|u|^{p-2}u\phi
dv_g-k\int_{B(x_i,\varepsilon)}|u|^{p-2}u\phi
d\sigma_g\\&&+\int_{B(x_i,\varepsilon)}f\phi
d\nu+\int_{B(x_i,\varepsilon)}K\phi d\pi
\end{eqnarray*}
this implies, after letting $\varepsilon$ tends to zero, that for
each $i\in J$, depending on wether $x_i\in M$ or $x_i\in\partial M$,
\begin{equation*}
    \mu_i\le f(x_i)\nu_i \text{ or }\mu_i\le K(x_i)\pi_i.
\end{equation*}
In particular, we assume that either $f(x_i)>0$ or $K(x_i)>0$,
because otherwise we get
$\mu_i=\nu_i=\pi_i=0$ and we are done.\\
Now, suppose that there exists $i\in I$ such that $\nu_i\neq0$ or
$\pi_i\neq0$ (depending on wether $x_i$ is in $M$ or $\partial M$),
then we get
\begin{eqnarray*}
\underset{(q,r\rightarrow(p^*,p^{**})} \lim\mu_{q,r}&=&
\underset{(q,r\rightarrow(p^*,p^{**})}\lim
\left((1-\frac{p}{q})\int_Mfv_{q,r}^qdv_g+
(1-\frac{p}{r})\int_{\partial M}Kv_{q,r}^rd\sigma_g\right) \\
&=&\frac{p}{n}\int_Mfu^{p^*}dv_g+\frac{p-1}{n-1}\int_{\partial
M}Ku^{p^{**}}d\sigma_g+\frac{p}{n}f(x_i)\nu_i+\frac{p-1}{n-1}K(x_i)\pi_i\\&\geq&
E(u)+\frac{p}{n}\left[\sup_Mf\right]^{1-\frac{n}{p}}K_1^{-n},
\end{eqnarray*}
which contradicts the hypothesis of the lemma. Thus, we get
$\mu_i=\nu_i=\pi_i=0$ and the sequence $v_{q,r}$ converges strongly
to $u$ in $L_{p^*}(M)$ and $L_{p^{**}}(\partial M)$.
\end{proof}
In the following, we give a sufficient condition in order to get
satisfied condition \eqref{eqt6.4} of lemma \ref{lem6.2}.
\begin{lem}\label{lem6.3}
Suppose that there exists a positive function $\Phi\in H^p_1(M)$such
that $I(\Phi)>0, \int_Mf\Phi^{p^*}dv_g>0, \int_{\partial
M}K\Phi^{p^{**}}d\sigma_g>0$ and
\begin{eqnarray*}
 0<\underset{\lambda\in[0,1]}\sup E(\lambda\Phi) &<&
E(u)+\frac{p}{n}\left[\sup_Mf\right]^{1-\frac{n}{p}}K_1^{-n}
\end{eqnarray*}
then  condition \eqref{eqt6.4} of lemma \ref{lem6.2} is satisfied.
\end{lem}
\begin{proof} Let  $u_{\ell_1},u_{\ell_2}$ be such that
$E_{q,r}(u_{\ell_1})=E_{q,r}(u_{\ell_2})=0$ and
$\ell_1=\|u_{\ell_1}\|^q_{q, M}+\|u_{\ell_1}\|^r_{r,\partial
M}<\ell_o<\ell_2=\|u_{\ell_2}\|^q_{q,
M}+\|u_{\ell_2}\|^r_{r,\partial M}$, where $\ell_o$ is such that
$\mu_{\ell_o}=\sup_{\ell\in]\ell_1,\ell_2[}\mu_{q,r,\ell}$. Then,
the curve $\lambda\rightarrow
E_{q,r}(\frac{\alpha\lambda-1}{\alpha-1} u_{\ell_2}), \alpha>1,
\lambda\in\left[\frac{1}{\alpha},\infty\right[$, starts from zero,
increases towards a positive maximum for $\lambda<1$ and
then decreases to minus infinity.\\
Now, suppose that there exists a function $\Phi\in H^p_1(M)$ such
that $\|\Phi\|^q_{q, M}+\|\Phi\|^r_{r,\partial M}>\ell_2$. Let
$\delta_o>\delta_1>1$ be two constants and consider the the curve
$\lambda\rightarrow E_{q,r}(
(1-\delta_o\lambda)(\delta_1\lambda-1)\Phi),
\lambda\in[\frac{1}{\delta_1},\infty]$, then we have,
\begin{eqnarray*}
\frac{dE_{q,r}((1-\delta_o\lambda)(\delta_1\lambda-1)\Phi)}{d\lambda}=
p(\delta_1+\delta_o)\left[|1-\delta_o\lambda||\delta_1\lambda-1|\right]^{p-1}
(1-\frac{2\delta_o\delta_1}{\delta_1+\delta_o}\lambda)\\
\left[\|\nabla_g\Phi\|^p_{p,M}+a\|\Phi\|^p_{p,M}+
k\|\Phi\|^p_{p,\partial
M}-\left(|1-\delta_o\lambda||\delta_1\lambda-1|\right)^{q-p}\int_Mf\Phi^qdv_g\right.\\
\left.-
\left(|1-\delta_o\lambda||\delta_1\lambda-1|\right)^{r-p}\int_{\partial
M}K\Phi^rdv_g\right].&
\end{eqnarray*}
Note $F_\Phi(\lambda)$ the function
\begin{eqnarray*}
F_\Phi(\lambda)=I(\phi)-
\left((1-\delta_o\lambda)(\delta_1\lambda-1)\right)^{q-p}\int_Mf\Phi^qdv_g\\
-\left((1-\delta_o\lambda)(\delta_1\lambda-1)\right)^{r-p}\int_{\partial
M}K\Phi^rd\sigma_g
\end{eqnarray*}
with
$\lambda\in[\frac{1}{\delta_o},\frac{1}{\delta_1}],I(\phi)>0,\int_Mf\Phi^qdv_g>0$
and $ \int_{\partial M}K\Phi^rd\sigma_g>0$.\\
Then, in the interval $(\frac{1}{\delta_o},\frac{1}{\delta_1})$,
there exist at most two values $\frac{1}{\delta_2}>
\frac{\delta_o+\delta_1}{2\delta_o\delta_1}>\frac{1}{\delta_3}$ of
$\lambda$ such that
$F_\Phi(\frac{1}{\delta_2})=F_\Phi(\frac{1}{\delta_2})=0,
F_\Phi(\lambda)>0,\lambda\in(\frac{1}{\delta_o},\frac{1}{\delta_3})\cup
(\frac{1}{\delta_2},\frac{1}{\delta_1})$ and
$F_\Phi(\lambda)<0,\lambda\in (\frac{1}{\delta_3},
\frac{1}{\delta_2}) $.\\
Hence, there exists $\delta_o>\delta_o'>\delta_1$ such that the
curve $\lambda\rightarrow E_{q,r}((1-\delta_o\lambda)
(\delta_1\lambda-1)\Phi)$ is positive for
$\lambda\in[\frac{1}{\delta_o'},\frac{1}{\delta_1}]$
and attains positive maximum at $\lambda_o=\frac{1}{\delta_2}$.\\
On the other hand, the curve $\lambda\rightarrow
E_{q,r}((1-\delta_o\lambda)u_{\ell_1})$, starts from zero, decreases
to  negative minimum in the interval $(0,\frac{1}{\delta_o})$ and
then increases to infinity.\\
 Now, take $\delta_o$ close to $\delta_o'$, $\delta_1$ as close to $1$
as \begin{equation*} E_{q,r}(\frac{\delta_1\lambda-1}{\delta_1-1}
u_{\ell_1})<\sup_{\lambda\in(\frac{1}{\delta'_o},\frac{1}{\delta_1})}E_{q,r}(
(1-\delta_o\lambda)(\delta_1\lambda-1)\Phi)=E_{q,r}(
\frac{1}{\delta_2}\Phi)
 \end{equation*} and
consider the curve
\begin{equation*}
g(\lambda)=\left\{
\begin{array}{ll}
(1-\delta_o\lambda)u_{\ell_1}, & 0\leq\lambda\leq\frac{1}{\delta_o} \\
(1-\delta_o\lambda)(\delta_1\lambda-1)\Phi,
&\frac{1}{\delta_o}\leq\lambda\leq\frac{1}{\delta_1}\\
\frac{\delta_1\lambda-1}{\delta_1-1} u_{\ell_2},&
\frac{1}{\delta_1}\leq\lambda\leq1.
                   \end{array}
                 \right.
\end{equation*}
Suppose that the condition of the lemma is satisfied, then in a
neighborhood  $V_{(p^{*},p^{**})}$of $(p^{*},p^{**})$ we can assume
for every $(q,r)\in V_{(p^{*},p^{**})}$ that
\begin{eqnarray*}
 0<\underset{\lambda\in[0,1]}\sup E_{q,r}(\lambda\Phi) &<&
E_{q,r}(u)+\frac{p}{n}\left[\sup_Mf\right]^{1-\frac{n}{p}}K_1^{-n}
\end{eqnarray*}
then
\begin{eqnarray*}
\mu_{q,r}&\leq&\underset{\lambda\in[o,1]}\sup
E_{q,r}(g(\lambda))\leq E_{q,r}(\frac{1}{4\delta_o}\Phi)\\&\leq&
E_{q,r}(u)+\frac{p}{n}\left[\sup_Mf\right]^{1-\frac{n}{p}}K_1^{-n}
\end{eqnarray*}
\end{proof}
As a result, we obtain the following lemma
\begin{lem}\label{lem6.4}
Suppose that the functions $f$ and $K$ satisfy
\begin{equation*}
\frac{p(n-1)}{n(p-1)}(\sup_Mf)^{-\frac{n-p}{p}}(\sup_{\partial
M}K)^\frac{n-p}{p-1}K_1^{-n}K_2^{\frac{p(n-1)}{p-1}}\le1
\end{equation*}
 and that there exists a function $\Phi\in H^p_1(M)$ with
 $I(\Phi)>0$, $\int_Mf\Phi^{p^*}dv_g>0$, $\int_{\partial
M}K\Phi^{p^{**}}d\sigma_g>0$ such that
\begin{eqnarray*}
  0<\underset{\lambda\in[0,1]}\sup E(\lambda\Phi) &<&
E(u)+\frac{p}{n}\left[\sup_Mf\right]^{1-\frac{n}{p}}K_1^{-n}
\end{eqnarray*}
Then none of the following cases
$\underset{(q,r)\rightarrow(p^*,p^{**})}\lim v_{q,r}=0$ or
\begin{equation*}
\lim \|v_{q,r}-u\|_{p,M}^p=\underset{(q,r)\rightarrow(p^*p^{**})}
\lim\|v_{q,r}-u\|_{p,\partial M}^p=0
\end{equation*}
 can occur.
\end{lem}
\begin{proof} If the sequence $v_{q,r}$ converges to zero function as $(q,r)$ goes to
$(p^*,p^{**})$, we can repeat the proof of lemma \ref{lem6.3}
without the term $E(u)$ we get
\begin{equation*}
\lim_{(q,r)\to(p^*,p^{**}}\mu_{q,r}\ge
\frac{p}{n}\left[\sup_Mf\right]^{1-\frac{n}{p}}K_1^{-n}
\end{equation*}
Since we have $E(u)<0$, under the hypothesis of the lemma  we get by
lemmas \ref{lem6.3} a contradiction, that is, $v_{q,r}$ does not
converge to zero function as $(q,r)$ goes to $(p^*,p^{**})$ and at
the same time the sequence $v_{q,r}$ can not satisfy that
\begin{equation*} \lim
\|v_{q,r}-u\|_{p,M}^p=\underset{(q,r)\rightarrow(p^*p^{**})}
\lim\|v_{q,r}-u\|_{p,\partial M}^p=0
\end{equation*}  since
$E_{q,r}(v_{q,r})=\lim_{(q,r)\to(p^*,p^{**}}\mu_{q,r}>0$.
\end{proof}
 Now the proof of theorem \ref{thm3.3} follows
\begin{proof}[\textbf{Proof of theorem \ref{thm3.3}}] Let $u$ be
the solution of problem \eqref{eq1.1} given by theorem\ref{thm3.1}.
Let $v_{q,r}$ the sequence of solution of the subcritical problem
\eqref{eq3.1} obtained by the mountain pass lemma. By lemma
\ref{lem6.1} it is bounded in $H^p_1(M)$, then after passing to a
subsequence we assume that $v_{q,r}$ converges , when $(q,r)$ goes
to $(p^*,p^{**})$ to a function $v$ weakly in $H^p_1(M)$, strongly
in $L_p(M)$ and $L_p(\partial M)$, almost everywhere in $M$ and in
the sense of trace on $\partial M$ and it converges to $v^{p^*-1}$
weakly in $L_{\frac{p^*}{p^*-1}}(M)$ and to $v^{p^{**}-1}$ weakly in
$L_{\frac{p^{**}}{p^{**}-1}}(\partial M)$.
Then we get that $v$ is a weak solution of problem \eqref{eq1.1}.\\
Under hypothesis of the theorem  we get by lemmas \ref{lem6.1},
\ref{lem6.2}, \ref{lem6.3} and \ref{lem6.4}  that $v\neq0$ and
$v\neq u$, that is problem \eqref{eq1.1} admits a second weak
solution.
\end{proof}
\bibliographystyle{amsplain}

\begin{thebibliography}{11}
\bibitem{Abdellaoui-Peral}B.Abdellaoui and I.Peral, On quasilinear
elliptic equations related to some Caffarelli-Khon-Nirenberg
inequalities, Communications on  Pure and Applied Analysis, Vol 2
N°4 (2003) 539-566.
\bibitem{Adam}R.A.Adams and J.J.F. Fournier, Sobolev spaces, Academic Press, Second
edition 2003 .
\bibitem{Ambrosetti-Li-MAlchiodi}, A. Ambrosetti, Y.Y.Li and
A.Malchiodi, Scalar curvature under boundary conditions, C.R.A.S.
Paris,t.330,Série 1,(2000),1013-1018.
\bibitem{Ambrosetti-Li-MAlchiodi1}A. Ambrosetti, Y.Y.Li and
A.Malchiodi, Yamabe and scalar curvature problems under boudary
conditions, Math.Ann.322(2002)667-699.
\bibitem{Aubin }T.Aubin, Some nonlinear problems in Riemannian
geometry, Springer,1998.
\bibitem{aubin-cotsiolis}T.Aubin and A. Cotsiolis, Equations non
linéaires avec le $p-$laplacien et la fonction exponentielle sur les
variétés riemanniennnes compactes, Bull.Sci.Math.124(2000)1-19.
\bibitem{Benalili} M.Benalili, Existence and multiplicity of solutions to elliptic
 equations of fourth order on compact manifolds, Dynamics of PDE, Vol(6)N°3(2009)203-225,
\bibitem{BenaliliMaliki1}M. Benalili and Y. Maliki, Generalized scalar curvature
type equation on complete Riemannian manifolds,
Elec.J.Diff.Eqns.146(2004)1-18.
\bibitem{BenaliliMaliki2} M. Benalili and Y. Maliki, Solving p-Laplacian on
complete Riemannian manifolds, Elec.J.Diff.Eqns.N°155(2006)1-9.
\bibitem{BenaliliMaliki3}M. Benalili and Y. Maliki, Generalized prescribed scalar
curvature type equation on compact manifolds of negative scalar
curvature, Rocky Mountain Journal of Mathematics Vol(5) N°37 (2007),
1399-1413.
\bibitem{Biezuner} R.J.Biezuner, Best constants in Sobolev trace
inequalities, Nonlinear Analysis 54 (2003) 575 – 589.
\bibitem{Biezuner1} R.J.Biezuner, Best constants, optimal Sobolev
inequalities on Riemannian manifolds and applications, Thesis, Univ.
New Jersy, 2003.
\bibitem{Bresis-Lieb} H. Bresis and E.H. Lieb, A relation between
pointwise convergence of functions and functionals, Proc. AMS.
88(1983)486-490.
\bibitem{cherrier} P.Cherrier, Problèmes de Neumann non linéaires sur
les variétés Riemanniennes, J.Func.An.57 N°2(1984)154-206.
\bibitem{cherrier1} P.Cherrier, Cas d'exception du théorème
d'inclusion de Sobolev sur les variétés Riemanniennes et
applications, Bull.Sci.Math.105(1981) 235-288.
\bibitem{cherrier2}P.Cherrier, Meilleurs constantes dans des
inégalités relatives aux espaces de Sobolev,
Bull.Sci.Math.108(1984)225-262.
\bibitem{cotsiolis-labropoulos} A.Cotsiolis and N.Labropoulos, A
Neumann problem with the $q-$laplacian on a solid torus in the
critical of supercritical case, Elec.J.Diff.Eqns,164(2007)1-18.
\bibitem{cotsiolis-labropoulos1} A.Cotsiolis and N.Labropoulos,
Best constants in Sobolev inequalities on manifolds with boundary in
the presence of symmetries and applications, Bull.Sci.Math. 132
(2008) 562–574.
\bibitem{Djadli-Malchiodi-Ahmedou}Z.Dajdli, A.Malchiodi and
M.O.Ahmedou, The prescribed boundary mean curvature problem on
$\mathbb{B}^4$, J.Diff.Eqns.206(2004)373-398.
\bibitem{Druet}O.Druet, Generalized scalar curvature type equations
on compact Riemannian manifolds, Proceedings of the Royal Society of
Edinburgh,130A(2000) 269-289.
\bibitem{Escobar1} J.Escobar, Conformal deformation of a Riemannian
metric to scalr flat metric with cosnatnt mean curvature on the
boundary. Ann.Math.136(1992)1-50 and 139(1994)749-750.
\bibitem{Escobar2} J.Escobar, The Yamabe probelm on manifolds with
boundary.J.Diff.Geo.35(1992)21-84.
\bibitem{Yasov-Runst} Y.IL’Yasov  and T.Runst,
Positive solutions for indefinite inhomgenous Neumann elliptic
probems, Electronic Journal of Differential Equations, N°57(2003)
1–21.
\bibitem{Lieberman}G.M.Lieberman, Boundary regularity for solutions of degenerate elliptic
equations, Nonlinear Analysis T.M.A. 12(1988) 1203–1219.
\bibitem{phozaev}S.I.Pohozaev, On an approach to Nonlinear equations. Doklady Acad. Sci. USSR
247 (1979), 1327-1331.
\bibitem{rauzy}A. Rauzy, Courbures scalaires des variétés d'invariant
conforme négatif, Trans.AMS.374 N°12(1995)4729-4745.
\end{thebibliography}

\end{document}